\documentclass [11pt,a4paper] {article}

\usepackage[cp1252]{inputenc}

\usepackage{amssymb}
\usepackage{amsmath}
\usepackage{amsfonts,amssymb}
\usepackage[dvips]{graphicx}
\usepackage{bbm}
\usepackage{enumerate}

\usepackage{amsthm}

\usepackage{cancel}

\DeclareMathAlphabet{\mathpzc}{OT1}{pzc}{m}{it}

\def\SmallColSep{\setlength{\arraycolsep}{1pt}}

\begin{document}

\title{The assumption of the Hilbert lattice in the case of\\a two-dimensional system}

\author{Arkady Bolotin\footnote{$Email: arkadyv@bgu.ac.il$\vspace{5pt}} \\ \textit{Ben-Gurion University of the Negev, Beersheba (Israel)}}

\maketitle

\begin{abstract}\noindent As it is known, the set of all closed linear subspaces of a Hilbert space together with a binary relation over the set represents the logic of the quantum propositions. It is also known that the lattices of the closed linear subspaces on a Hilbert space of dimension 3 or greater do not have a prime filter, hence those lattices do not allow a valuation map. In contrast to that, for qubits it is easy to find prime filters in the Hilbert lattice.\\

\noindent This begs the question: What assumption(s) related to the lattices of the closed linear subspaces should be added or altered to preclude the bivaluation map in the two-dimensional case? The presented paper offers the answer to this question.
\\

\end{abstract}

\section{Introduction and preliminaries}  

\noindent Consider \textit{the complete lattice $(\mathcal{L}(\mathbb{C}^2),\le)$}, where $\le$ denotes the partial ordering over $\mathcal{L}(\mathbb{C}^2)$, the set of all \textit{the closed linear subspaces} of the two-dimension Hilbert space $\mathcal{H} = \mathbb{C}^2$, namely,\smallskip

\begin{equation}  
   \mathcal{L}(\mathbb{C}^2)
   =
   \Bigg\{
      \{0\}
      ,
      \left\{
         \!\left[
            \begingroup\SmallColSep
            \begin{array}{r}
               a \\
               a
            \end{array}
            \endgroup
         \right]\!
      \right\}
      ,
      \left\{
         \!\left[
            \begingroup\SmallColSep
            \begin{array}{r}
               a \\
              -a
            \end{array}
            \endgroup
         \right]\!
      \right\}
      ,
      \left\{
         \!\left[
            \begingroup\SmallColSep
            \begin{array}{r}
               ia \\
               a
            \end{array}
            \endgroup
         \right]\!
      \right\}
      ,
      \left\{
         \!\left[
            \begingroup\SmallColSep
            \begin{array}{r}
               a \\
              ia
            \end{array}
            \endgroup
         \right]\!
      \right\}
      ,
      \left\{
         \!\left[
            \begingroup\SmallColSep
            \begin{array}{r}
               a \\
               0
            \end{array}
            \endgroup
         \right]\!
      \right\}
      ,
      \left\{
         \!\left[
            \begingroup\SmallColSep
            \begin{array}{r}
               0 \\
               a
            \end{array}
            \endgroup
         \right]\!
      \right\}
      ,
      \mathbb{C}^2
   \Bigg\}
   \;\;\;\;  ,
\end{equation}
\smallskip

\noindent in which $a\in \mathbb{R}$. The said subspaces are \textit{the column spaces} (or \textit{images}, or \textit{ranges}) of the projection operators $\hat{P}^{(Q)}_n$ on $\mathbb{C}^2$, i.e.,\smallskip

\begin{equation}  
   \left\{
      \!\left[
         \begingroup\SmallColSep
         \begin{array}{r}
            \cdot \\
            \cdot
         \end{array}
         \endgroup
      \right]\!
   \right\}
   =
   \mathrm{ran}(\hat{P}^{(Q)}_n)
   \;\;\;\;  ,
\end{equation}
\smallskip

\noindent where $\hat{P}^{(Q)}_n$ are defined by the formula

\vspace*{1mm}

\begin{equation}  
   \hat{P}^{(Q)}_n
   =
   \frac{1}{2}
      \left[
         \begin{array}{l l}
            1+(-1)^n (\delta_{Q0}-\delta_{Q3})
            &
            (-1)^n (-\delta_{Q1}+i\delta_{Q2})
            \\
            (-1)^n (-\delta_{Q1}-i\delta_{Q2})
            &
            1+(-1)^n (\delta_{Q0}+\delta_{Q3})
         \end{array}
      \right]
   \;\;\;\;  ,
\end{equation}
\smallskip

\noindent in which $\delta_{ab}$ is the Kronecker delta, $n \in \{1,2\}$ and $Q=\{0,1,2,3\}$. According to this formula,\\

\begin{equation}  
   \{0\}
   =
   \mathrm{ran}(\hat{P}^{(0)}_1)
   =
   \mathrm{ran}(\hat{0})
   \;\;\;\;  ,
\end{equation}

\begin{equation}  
   \mathbb{C}^2
   =
   \mathrm{ran}(\hat{P}^{(0)}_2)
   =
   \mathrm{ran}(\hat{1})
   \;\;\;\;  ,
\end{equation}
\smallskip

\noindent where $\hat{0}$ and $\hat{1}$ are the zero and identity projection operators respectively.\\

\noindent The partial ordering $\le$ on $\mathcal{L}(\mathbb{C}^2)$ is defined by\smallskip

\begin{equation} \label{A1} 
   \mathrm{ran}(\hat{0})
   \le
   \mathrm{ran}(\hat{P}^{(Q)}_n)
   \iff
   \{0\}
   \subseteq
   \mathrm{ran}(\hat{P}^{(Q)}_n)
   \;\;\;\;  ,
\end{equation}

\begin{equation} \label{A2} 
   \mathrm{ran}(\hat{P}^{(Q)}_n)
   \le
   \mathrm{ran}(\hat{1})
   \iff
   \mathrm{ran}(\hat{P}^{(Q)}_n)
   \subseteq
   \mathbb{C}^2
   \;\;\;\;  .
\end{equation}
\smallskip

\noindent As stated by the definition of a complete lattice \cite{Nation, Davey}, each two-element subset of $\mathcal{L}(\mathbb{C}^2)$, namely, $\{ \mathrm{ran}(\hat{P}^{(Q)}_n), \mathrm{ran}(\hat{P}^{(R)}_m) \}\subset \mathcal{L}(\mathbb{C}^2)$, where $m \in \{1,2\}$, $R \in \{0,1,2,3\}$ such that $\mathrm{ran}(\hat{P}^{(Q)}_n) \neq \mathrm{ran}(\hat{P}^{(R)}_m)$, has \textit{a meet} denoted by $\mathrm{ran}(\hat{P}^{(Q)}_n) \wedge \mathrm{ran}(\hat{P}^{(R)}_m)$ and \textit{a join} denoted by $\mathrm{ran}(\hat{P}^{(Q)}_n) \vee \mathrm{ran}(\hat{P}^{(R)}_m)$. In symbols,\smallskip

\begin{equation} \label{B} 
   \mathrm{ran}(\hat{P}^{(Q)}_n)
   \wedge
   \mathrm{ran}(\hat{P}^{(R)}_m)
   =
   \mathrm{ran}(\hat{P}^{(Q)}_n)
   \cap
   \mathrm{ran}(\hat{P}^{(R)}_m)
   \in
   \mathcal{L}(\mathbb{C}^2)
   \;\;\;\;  ,
\end{equation}

\begin{equation} \label{C} 
   \mathrm{ran}(\hat{P}^{(Q)}_n)
   \vee
   \mathrm{ran}(\hat{P}^{(R)}_m)
   =
   \left(
   \mathrm{ran}(\hat{P}^{(Q)}_{n^\prime \neq n})
   \cap
   \mathrm{ran}(\hat{P}^{(R)}_{m^\prime \neq m})
   \right)^{\perp}
   \in
   \mathcal{L}(\mathbb{C}^2)
   \;\;\;\;  ,
\end{equation}
\smallskip

\noindent where $(\cdot)^{\perp}$ stands for the orthogonal complement of $(\cdot)$ in a way that $(\mathrm{ran}(\hat{P}^{(Q)}_{n}))^{\perp}\! = \mathrm{ran}(\hat{P}^{(Q)}_{n^\prime \neq n}) \in \mathcal{L}(\mathbb{C}^2)$. In view of (\ref{A1}) and (\ref{A2}), it follows then\smallskip

\begin{equation}  
   \mathrm{ran}(\hat{0})
   \le
   \mathrm{ran}(\hat{P}^{(Q)}_n)
   \iff
   \mathrm{ran}(\hat{0})
   \wedge
   \mathrm{ran}(\hat{P}^{(Q)}_n)
   =
   \mathrm{ran}(\hat{0})
   \;\;\;\;  ,
\end{equation}

\begin{equation}  
   \mathrm{ran}(\hat{P}^{(Q)}_n)
   \le
   \mathrm{ran}(\hat{1})
   \iff
   \mathrm{ran}(\hat{P}^{(Q)}_n)
   \wedge
   \mathrm{ran}(\hat{1})
   =
   \mathrm{ran}(\hat{P}^{(Q)}_n)
   \;\;\;\;  .
\end{equation}
\smallskip

\noindent Consider $\mathcal{F}(\mathrm{ran}(\hat{P}^{(W)}_k))$, a subset in the lattice $(\mathcal{L}(\mathbb{C}^2),\le)$ defined by\smallskip

\begin{equation}  
   \mathcal{F}\left(\mathrm{ran}(\hat{P}^{(W)}_k)\right)
   =
   \mathcal{L}(\mathbb{C}^2)
   \setminus
   \left\{
   \mathrm{ran}(\hat{P}^{(W)}_k)
   \right\}
   \;\;\;\;  ,
\end{equation}
\smallskip

\noindent where $k \in \{1,2\}$ and $W \in \{1,2,3\}$, in other words, $\mathrm{ran}(\hat{P}^{(W)}_k)$ is the nontrivial closed subspace belonging to $\mathcal{L}(\mathbb{C}^2)$. According to (\ref{B}), $\mathcal{F}(\mathrm{ran}(\hat{P}^{(W)}_k))$ is downward directed, that is, for all $\mathrm{ran}(\hat{P}^{(Q)}_n)$ and $\mathrm{ran}(\hat{P}^{(R)}_m)$ in $\mathcal{F}(\mathrm{ran}(\hat{P}^{(W)}_k))$ one finds\smallskip

\begin{equation}  
   \mathrm{ran}(\hat{P}^{(Q)}_n)
   \wedge
   \mathrm{ran}(\hat{P}^{(R)}_m)
   \in
   \mathcal{F}\left(\mathrm{ran}(\hat{P}^{(W)}_k)\right)
   \;\;\;\;  ,
\end{equation}
\smallskip

\noindent which means that the subset $\mathcal{F}(\mathrm{ran}(\hat{P}^{(W)}_k))$ is \textit{a filter} in the lattice $(\mathcal{L}(\mathbb{C}^2),\le)$ \cite{Burris}. Moreover, in accordance with (\ref{C}), for all $\mathrm{ran}(\hat{P}^{(Q)}_n),\, \mathrm{ran}(\hat{P}^{(W)}_k) \in \mathcal{L}(\mathbb{C}^2)$ such that\smallskip

\begin{equation}  
   \mathrm{ran}(\hat{P}^{(Q)}_n)
   \vee
   \mathrm{ran}(\hat{P}^{(W)}_k)
   \in
   \mathcal{F}\left(\mathrm{ran}(\hat{P}^{(W)}_k)\right)
   \;\;\;\;   
\end{equation}
\smallskip

\noindent one has $\mathrm{ran}(\hat{P}^{(Q)}_n) \in \mathcal{F}(\mathrm{ran}(\hat{P}^{(W)}_k))$, which means that $\mathcal{F}(\mathrm{ran}(\hat{P}^{(W)}_k))$ is also \textit{a prime filter} in $(\mathcal{L}(\mathbb{C}^2),\le)$, whereas\smallskip

\begin{equation}  
   \mathcal{I}\left(\mathrm{ran}(\hat{P}^{(W)}_k)\right)
   =
   \mathcal{L}(\mathbb{C}^2)
   \setminus
   \mathcal{F}\left(\mathrm{ran}(\hat{P}^{(W)}_k)\right)
   =
   \left\{
   \mathrm{ran}(\hat{P}^{(W)}_k)
   \right\}
   \;\;\;\;   
\end{equation}
\smallskip

\noindent is \textit{a prime ideal} in $(\mathcal{L}(\mathbb{C}^2),\le)$ \cite{Burris}.\\

\noindent Because the lattice $(\mathcal{L}(\mathbb{C}^2),\le)$ contains the prime filters, one can define \textit{the lattice homomorphism $h\!:\,(\mathcal{L}(\mathbb{C}^2),\le) \mapsto (\{ \{0\},\mathbb{C}^2 \}, \le)$} by\smallskip

\begin{equation}  
   h\left(\mathrm{ran}(\hat{P}^{(Q)}_n)\right)
   =
   \left\{
      \begin{array}{l}
         \,\,\mathbb{C}^2, \;\;\;\;  \mathrm{if} \;\mathrm{ran}(\hat{P}^{(Q)}_n) \in \mathcal{I}\left(\mathrm{ran}(\hat{P}^{(W)}_k)\right) \\
         \{0\},                       \;\;\;\;  \mathrm{if} \;\mathrm{ran}(\hat{P}^{(Q)}_n) \in \mathcal{F}\left(\mathrm{ran}(\hat{P}^{(W)}_k)\right)
      \end{array}
   \right.
   \;\;\;\;   
\end{equation}
\smallskip

\noindent and, consequently, the bivaluation map $v\!:\mathcal{L}(\mathbb{C}^2) \mapsto \{0,1\}$ by\smallskip

\begin{equation} \label{D} 
   v\left(\mathrm{ran}(\hat{P}^{(Q)}_n)\right)
   =
   \left\{
      \begin{array}{l}
         1, \;\;\;\;  \mathrm{if} \;h\left(\mathrm{ran}(\hat{P}^{(Q)}_n)\right) = \mathbb{C}^2 \\
         0, \;\;\;\;  \mathrm{if} \;h\left(\mathrm{ran}(\hat{P}^{(Q)}_n)\right) = \{0\}
      \end{array}
   \right.
   \;\;\;\;  ,
\end{equation}
\smallskip

\noindent where 1 and 0 represent the truth and falsity, respectively.\\

\noindent Given that any $\hat{P}^{(Q)}_n$ is the identity operator on $\mathrm{ran}(\hat{P}^{(Q)}_n)$ and the zero operator on $\mathrm{ker}(\hat{P}^{(Q)}_n) = \mathrm{ran}(\neg\hat{P}^{(Q)}_n)$, where $\neg\hat{P}^{(Q)}_n = \hat{1} - \hat{P}^{(Q)}_n$, the valuation (\ref{D}) can be modified further as\smallskip

\begin{equation} \label{E} 
   v\left(\hat{P}^{(Q)}_n \right)
   =
   \left\{
      \begin{array}{l}
         1, \;\;\;\;  |\Psi\rangle \in \mathrm{ran}(\hat{P}^{(Q)}_n)\\
         0, \;\;\;\;  |\Psi\rangle \in \mathrm{ran}(\neg\hat{P}^{(Q)}_n)
      \end{array}
   \right.
   \;\;\;\;  ,
\end{equation}
\smallskip

\noindent where $|\Psi\rangle$ describes the state in which the two-dimensional system is prepared. Due to (\ref{E}), all the propositions corresponding to the qubit projection operators $\hat{P}^{(Q)}_n$ must obey \textit{the principle of bivalence} (according to which every proposition can be either true or false \cite{Beziau}).\\

\noindent On the other hand, it is known that \textit{a Hilbert lattice} – i.e., a lattice of all the closed linear subspaces of a Hilbert space $\mathcal{H}$ – does not have a prime filter if $\dim(\mathcal{H}) \ge 3$. As a result, there does not exist a bivaluation map $v\!:\mathcal{L}(\mathcal{H}) \mapsto \{0,1\}$ on the Hilbert lattice $(\mathcal{L}(\mathcal{H}), \le)$ with $\dim(\mathcal{H}) \ge 3$ \cite{Redei}.\\

\noindent This means that if two-dimensional systems are not excluded from the domain of validity of quantum mechanics, the valuation (\ref{D}) should be regarded as \textit{physically unsound}. Accordingly, the question is, what assumption(s) related to the lattices of the closed linear subspaces should be added or altered to preclude the bivaluation map $v\!:\mathcal{L}(\mathbb{C}^2) \mapsto \{0,1\}$?\\

\noindent The presented paper offers the answer to this question.\\

\section{Lattices of subspaces of the Hilbert space}  

\noindent Mathematically, the reason for the bivalence of all the propositions corresponding to the qubit projection operators $\hat{P}^{(Q)}_n$ is that \textit{the assumption of the Hilbert lattice is not strong enough} to exclude the prime filters $\mathcal{F}(\mathrm{ran}(\hat{P}^{(W)}_k))$ in the partially ordered set $\mathcal{L}(\mathbb{C}^2)$.\\

\noindent To see this, let us formally describe a set of the closed linear subspaces of the Hilbert space using set-builder notation:\smallskip

\begin{equation} \label{F1} 
   \left\{
      \mathrm{ran}(\hat{P}^{(Q)}_n)
      :\;
      \Phi\left( \mathrm{ran}(\hat{P}^{(Q)}_n) \right)
   \right\}
   \;\;\;\;  .
\end{equation}
\smallskip

\noindent If the logical predicate $\Phi$, i.e., the rule defining the set, holds for all the closed subspaces of $\mathbb{C}^2$, namely,\smallskip

\begin{equation} \label{F2} 
   \Phi\left( \mathrm{ran}(\hat{P}^{(Q)}_n) \right)
   =
   \mathrm{ran}(\hat{P}^{(Q)}_n)
   \subseteq
   \mathbb{C}^2
   \;\;\;\;  ,
\end{equation}
\smallskip

\noindent the definition (\ref{F1}) is read\smallskip

\begin{equation}  
   \left\{
      \mathrm{ran}(\hat{P}^{(Q)}_n)
      :\;
      \Phi\left( \mathrm{ran}(\hat{P}^{(Q)}_n) \right)
   \right\}
   =
   \mathcal{L}(\mathbb{C}^2)
   \;\;\;\;  .
\end{equation}
\smallskip

\noindent However, an assumption stronger than (\ref{F2}) is possible here. Particularly, the predicate $\Phi$ holds only for those subspaces $\mathrm{ran}(\hat{P}^{(Q)}_n) \subseteq \mathbb{C}^2$ that are \textit{invariant} under the nontrivial projection operator $\hat{P}^{(W)}_k$, i.e., $\hat{P}^{(W)}_k\!: \mathrm{ran}(\hat{P}^{(Q)}_n) \mapsto \mathrm{ran}(\hat{P}^{(Q)}_n)$. Explicitly,\smallskip

\begin{equation} \label{G} 
   \Phi\left( \mathrm{ran}(\hat{P}^{(Q)}_n) \right)
   =
   \left(
      \mathrm{ran}(\hat{P}^{(Q)}_n)
      \subseteq
      \mathbb{C}^2
   \right)
   \wedge
   \left(
      \hat{P}^{(W)}_k\!
      :\;
      \mathrm{ran}(\hat{P}^{(Q)}_n)
      \mapsto
      \mathrm{ran}(\hat{P}^{(Q)}_n)
   \right)
   \;\;\;\;  ,
\end{equation}
\smallskip

\noindent where\smallskip

\begin{equation}  
   \left(
      \hat{P}^{(W)}_k\!
      :\;
      \mathrm{ran}(\hat{P}^{(Q)}_n)
      \mapsto
      \mathrm{ran}(\hat{P}^{(Q)}_n)
   \right)
   \in
   \left\{
      |\Psi\rangle
      \in
      \mathrm{ran}(\hat{P}^{(Q)}_n)
      :\,
      \hat{P}^{(W)}_k |\Psi\rangle
      \subseteq
      \mathrm{ran}(\hat{P}^{(Q)}_n)      
   \right\}
   \;\;\;\;  .
\end{equation}
\smallskip

\noindent The set of the subspaces $\mathrm{ran}(\hat{P}^{(Q)}_n)$ that satisfy the rule (\ref{G}) is as follows:\smallskip

\begin{equation}  
   \left\{
      \mathrm{ran}(\hat{0})
      ,\,
      \mathrm{ran}(\hat{P}^{(W)\!}_{k})
      ,\,
      \mathrm{ran}(\neg\hat{P}^{(W)\!}_{k})
      ,\,
      \mathrm{ran}(\hat{1})\!
   \right\}
   \equiv
   \mathcal{L}^{(W)}_{k}
   \;\;\;\;  .
\end{equation}
\smallskip

\noindent To be sure, let $|\Psi\rangle \in \mathrm{ran}(\hat{P}^{(W)}_{k})$. As $\hat{P}^{(W)}_{k} |\Psi\rangle = |\Psi\rangle$, one gets $\hat{P}^{(W)}_{k} |\Psi\rangle \in \mathrm{ran}(\hat{P}^{(W)}_{k})$, and so $\hat{P}^{(W)}_{k}\!: \mathrm{ran}(\hat{P}^{(W)}_{k}) \mapsto  \mathrm{ran}(\hat{P}^{(W)}_{k})$. Dually, let $|\Psi\rangle \in \mathrm{ran}(\neg\hat{P}^{(W)}_{k})$. This means that $\hat{P}^{(W)}_{k} |\Psi\rangle = 0$. On the other hand, $0 \in \mathrm{ran}(\neg\hat{P}^{(W)}_{k})$, which implies $\hat{P}^{(W)}_{k}\!: \mathrm{ran}(\neg\hat{P}^{(W)}_{k}) \mapsto  \mathrm{ran}(\neg\hat{P}^{(W)}_{k})$. Furthermore, the subspace $\mathrm{ran}(\hat{1}) = \mathcal{H}$ as well as the zero subspace $\mathrm{ran}(\hat{0}) = \{0\}$ are \textit{the trivially invariant subspaces} for any projection operator.\\

\noindent Define invariant subspaces for the sets of the operators as subspaces invariant for each operator in the set. Then, the invariant subspaces for the sets of the commutable projectors $\Sigma^{(W)} = \{P^{(W)}_k\}_{k=1}^2$ are\smallskip

\begin{equation}  
   \mathcal{L}^{(W)}
   =
   \bigwedge_{k \in \{1,2\}}
      \!\!\!\!
      \mathcal{L}^{(W)}_{k}
   =
   \bigcap_{k \in \{1,2\}}
      \!\!\!\!
      \mathcal{L}^{(W)}_{k}
   \;\;\;\;  ,
\end{equation}
\smallskip

\noindent or, explicitly,\smallskip

\begin{equation}  
   \left\{
      \mathrm{ran}(\hat{0})
      ,\,
      \mathrm{ran}(\hat{P}^{(1)\!}_{1})
      ,\,
      \mathrm{ran}(\hat{P}^{(1)\!}_{2})
      ,\,
      \mathrm{ran}(\hat{1})\!
   \right\}
   =
   \Bigg\{
      \{0\}
      ,
      \left\{
         \!\left[
            \begingroup\SmallColSep
            \begin{array}{r}
               a \\
               a
            \end{array}
            \endgroup
         \right]\!
      \right\}
      ,
      \left\{
         \!\left[
            \begingroup\SmallColSep
            \begin{array}{r}
               a \\
              -a
            \end{array}
            \endgroup
         \right]\!
      \right\}
      ,
      \mathbb{C}^2
   \Bigg\}
   \equiv
   \mathcal{L}^{(1)}
   \;\;\;\;  ,
\end{equation}

\begin{equation}  
   \left\{
      \mathrm{ran}(\hat{0})
      ,\,
      \mathrm{ran}(\hat{P}^{(2)\!}_{1})
      ,\,
      \mathrm{ran}(\hat{P}^{(2)\!}_{2})
      ,\,
      \mathrm{ran}(\hat{1})\!
   \right\}
   =
   \Bigg\{
      \{0\}
      ,
      \left\{
         \!\left[
            \begingroup\SmallColSep
            \begin{array}{r}
              ia \\
               a
            \end{array}
            \endgroup
         \right]\!
      \right\}
      ,
      \left\{
         \!\left[
            \begingroup\SmallColSep
            \begin{array}{r}
               a \\
              ia
            \end{array}
            \endgroup
         \right]\!
      \right\}
      ,
      \mathbb{C}^2
   \Bigg\}
   \equiv
   \mathcal{L}^{(2)}
   \;\;\;\;  ,
\end{equation}

\begin{equation}  
   \left\{
      \mathrm{ran}(\hat{0})
      ,\,
      \mathrm{ran}(\hat{P}^{(3)\!}_{1})
      ,\,
      \mathrm{ran}(\hat{P}^{(3)\!}_{2})
      ,\,
      \mathrm{ran}(\hat{1})\!
   \right\}
   =
   \Bigg\{
      \{0\}
      ,
      \left\{
         \!\left[
            \begingroup\SmallColSep
            \begin{array}{r}
               a \\
               0
            \end{array}
            \endgroup
         \right]\!
      \right\}
      ,
      \left\{
         \!\left[
            \begingroup\SmallColSep
            \begin{array}{r}
               0 \\
               a
            \end{array}
            \endgroup
         \right]\!
      \right\}
      ,
      \mathbb{C}^2
   \Bigg\}
   \equiv
   \mathcal{L}^{(3)}
   \;\;\;\;  .
\end{equation}
\smallskip

\noindent The elements of $\mathcal{L}^{(W)}$ form the complete lattices $(\mathcal{L}^{(W)}, \le)$ called \textit{invariant-subspace lattices} \cite{Radjavi}. Certainly, each $(\mathcal{L}^{(W)}, \le)$ is a lattice since every $\mathcal{L}^{(W)}$ has the greatest element, $\mathrm{ran}(\hat{1})$, and every pair of the elements of $\mathcal{L}^{(W)}$ has the meet that is the element of $\mathcal{L}^{(W)}$. Besides, every finite lattice is complete.\\

\noindent There is no subset containing $\mathrm{ran}(\hat{P}^{(W^\prime \neq W)}_{k \in \{1,2\}})$ in the partially ordered set $\mathcal{L}^{(W)}$, as well as there is no subset containing $\mathrm{ran}(\hat{P}^{(W)}_{k \in \{1,2\}})$ in the partially ordered set $\mathcal{L}^{(W^\prime \neq W)}$. This can also be worded by saying that each set $\mathcal{L}^{(W)}$ contains only the closed subspaces belonging to mutually commutable projection operators.\\

\noindent Consider\smallskip

\begin{equation}  
   \mathcal{L}
   =
   \bigwedge_{W \in \{1,2,3\}}
      \!\!\!\!
      \mathcal{L}^{(W)}
   =
   \bigcap_{W \in \{1,2,3\}}
      \!\!\!\!
      \mathcal{L}^{(W)}
   \;\;\;\;  ,
\end{equation}
\smallskip

\noindent i.e., the set of the invariant subspaces \textit{invariant under each nontrivial projection operator} $\hat{P}^{(W)}_{k}$ on the Hilbert space $\mathbb{C}^2$. Since the collection of all the nontrivial projection operators on $\mathbb{C}^2$\smallskip

\begin{equation}  
   \Sigma
   =
   \left\{
      \Sigma^{(W)}
   \right\}^3_{W=1}
   =
   \left\{
      \left\{
         \hat{P}^{(W)}_{k}
      \right\}^2_{k=1}
   \right\}^3_{W=1}
   \;\;\;\;   
\end{equation}
\smallskip

\noindent spans $\mathbb{C}^2$ and, hence, equals $\mathcal{A}(\mathbb{C}^2)$, \textit{the algebra of all linear transformations on $\mathbb{C}^2$}, the set $\mathcal{L}$ must be \textit{irreducible} in accordance with Burnside's theorem on incommutable algebras \cite{Burnside, Rosenthal, Lomonosov, Shapiro}, that is,\smallskip

\begin{equation}  
   \Sigma
   =
   \mathcal{A}(\mathcal{H})
   \;
   \implies
   \;
   \mathcal{L}
   =
   \left\{
      \{0\}
      ,
      \mathbb{C}^2
   \right\}
   \;\;\;\;  .
\end{equation}
\smallskip

\noindent Consequently, there is no subset containing $\mathrm{ran}(\hat{P}^{(W)}_{k \in \{1,2\}})$ and $\mathrm{ran}(\hat{P}^{(W^\prime \neq W)}_{k \in \{1,2\}})$ in the lattice $(\mathcal{L}, \le)$. One can conclude then that the subspaces $\mathrm{ran}(\hat{P}^{(W)}_{k \in \{1,2\}})$ can belong only to $\mathcal{L}^{(W)}$ while the subspaces $\mathrm{ran}(\hat{P}^{(W^\prime \neq W)}_{k \in \{1,2\}})$ only to $\mathcal{L}^{(W^\prime \neq W)}$. As a result, the nontrivial elements of the posets $\mathcal{L}^{(W)}$ and $\mathcal{L}^{(W^\prime \neq W)}$ cannot meet each other. In symbols,\smallskip

\begin{equation}  
   \mathrm{ran}(\hat{P}^{(W)}_{k \in \{1,2\}})
   \in
   \mathcal{L}^{(W)}
   ,\,
   \mathrm{ran}(\hat{P}^{(W^\prime \neq W)}_{k \in \{1,2\}})
   \in
   \mathcal{L}^{(W^\prime \neq W)}
   \,
   \implies
   \,
   \mathrm{ran}(\hat{P}^{(W)}_{k \in \{1,2\}})
   \;\cancel{\;\wedge\;}\;
   \mathrm{ran}(\hat{P}^{(W^\prime \neq W)}_{k \in \{1,2\}})
   \;\;\;\;  ,
\end{equation}
\smallskip

\noindent where the cancelation of $\wedge$ indicates that the meet operation cannot be defined for the two elements $\mathrm{ran}(\hat{P}^{(W)}_{k \in \{1,2\}})$ and $\mathrm{ran}(\hat{P}^{(W^\prime \neq W)}_{k \in \{1,2\}})$ of the different lattices $(\mathcal{L}^{(W)}, \le)$ and $(\mathcal{L}^{(W^\prime \neq W)}, \le)$. This means the following: Suppose that the filter $\mathcal{F}(\mathrm{ran}(\hat{P}^{(W)}_k))$ exists in the lattice $(\mathcal{L}_W, \le)$ and thus $v(\mathrm{ran}(\hat{P}^{(W)}_n)) \in \{0,1\}$, i.e.,\smallskip

\begin{equation}  
   v\left(\mathrm{ran}(\hat{P}^{(W)}_n)\right)
   =
   \left\{
      \begin{array}{l}
         1,                       \;\;\;\;  \mathrm{if} \;\mathrm{ran}(\hat{P}^{(W)}_n) \in \mathcal{I}\left(\mathrm{ran}(\hat{P}^{(W)}_k)\right) \\
         0,                       \;\;\;\;  \mathrm{if} \;\mathrm{ran}(\hat{P}^{(W)}_n) \in \mathcal{F}\left(\mathrm{ran}(\hat{P}^{(W)}_k)\right)
      \end{array}
   \right.
   \;\;\;\;  .
\end{equation}
\smallskip

\noindent At the same time, $v(\mathrm{ran}(\hat{P}^{(W^\prime \neq W)}_m))$ is neither 0 nor 1 because the said filter cannot exist in the lattice $(\mathcal{L}^{(W^\prime \neq W)}, \le)$. In symbols,\smallskip

\begin{equation}  
   \mathrm{ran}(\hat{P}^{(W^\prime \neq W)}_m) \in \cancel{\:\;\mathcal{F}\!\left(\mathrm{ran}(\hat{P}^{(W)}_k)\right)\;\;}
   \;
   \implies
   \;
   v\left(\mathrm{ran}(\hat{P}^{(W^\prime \neq W)}_m)\right)
   \;\cancel{\!\in\!}\;
   \{0,1\}
   \;\;\;\;  ,
\end{equation}
\smallskip

\noindent where the cancelation of $\mathcal{F}(\mathrm{ran}(\hat{P}^{(W)}_k))$ indicates that it cannot be a subset of the partially ordered set $\mathcal{L}^{(W^\prime \neq W)}$.\\

\noindent Correspondingly, the bivaluation map\smallskip

\begin{equation}  
   v
   :
   \left\{
      \mathcal{L}^{(1)}
      ,
      \mathcal{L}^{(2)}
      ,
      \mathcal{L}^{(3)}
   \right\}
   \;
   \mapsto
   \;
   \{0,1\}
   \;\;\;\;   
\end{equation}
\smallskip

\noindent is precluded because if $v\!:\mathcal{L}^{(W)} \mapsto \{0,1\}$, then $v\!:\mathcal{L}^{(W^\prime \neq W)} \;\;\cancel{\!\!\mapsto\!\!}\;\; \{0,1\}$.\\

\section{Concluding remarks}  

\noindent Even though the assertion that the set of \textit{all} the closed linear subspaces of a Hilbert space form a lattice looks mathematically unassuming, it appears to be too presumptuous from the physical point of view.\\

\noindent Firstly, this assertion introduces the meets and the joins of the closed subspaces associated with the \textit{incommutable} projection operators, such as\smallskip

\begin{equation} \label{J} 
   \mathrm{ran}(\hat{P}^{(W)}_n)
   \wedge
   \mathrm{ran}(\hat{P}^{(W^\prime \neq W)}_m)
   =
   \{0\}
   \;\;\;\;  ,
\end{equation}

\begin{equation} \label{K} 
   \mathrm{ran}(\hat{P}^{(W)}_{k \neq n})
   \wedge
   \mathrm{ran}(\hat{P}^{(W^\prime \neq W)}_m)
   =
   \{0\}
   \;\;\;\;  ,
\end{equation}
\smallskip

\noindent where $\hat{P}^{(W)}_n$ as well as $\hat{P}^{(W)}_{k \neq n}$ does not commute with $\hat{P}^{(W^\prime \neq W)}_m$. Substituting (\ref{J})  and (\ref{K}) into the distributive axiom\smallskip

\begin{equation}  
   \left(
      \mathcal{K}
      \vee
      \mathcal{M}
   \right)
   \wedge
   \mathcal{O}
   =
   \left(
      \mathcal{K}
      \wedge
      \mathcal{O}
   \right)
   \vee
   \left(
      \mathcal{M}
      \wedge
      \mathcal{O}
   \right)
   \;\;\;\;  ,
\end{equation}
\smallskip

\noindent where $\mathcal{K} = \mathrm{ran}(\hat{P}^{(W)}_n)$, $\mathcal{M} = \mathrm{ran}(\hat{P}^{(W)}_{k \neq n})$, and $\mathcal{O} = \mathrm{ran}(\hat{P}^{(W^\prime \neq W)}_m)$, one immediately finds a contradiction, namely, $\mathrm{ran}(\hat{P}^{(W^\prime \neq W)}_m) = \{0\}$.\\

\noindent Thus, to maintain that the collection of all the closed subspaces $\mathcal{L}(\mathcal{H})$ form a lattice, it is necessary to give up \textit{distributivity}. However, removing the distributive laws of propositional logic from quantum mechanics is equivalent to the assumption that quantum mechanics requires no less than a revolution in our understanding of logic \textit{per se}.\\

\noindent Secondly (and more importantly), the assumption of the Hilbert lattice $(\mathcal{L}(\mathcal{H}), \le)$ brings about the anomalous and \textit{physically unjustifiable} status of the two-dimensional case.\\

\noindent Therefore, as it is argued in the presented paper, the assumption of the Hilbert lattice $(\mathcal{L}(\mathcal{H}), \le)$ should be replaced by the assumption of the invariant-subspace lattices $(\mathcal{L}^{(W)}, \le)$. Despite its being stronger mathematically, the latter assumption does not bring in \textit{new physical hypotheses}. According to it, all the closed linear subspaces of $\mathbb{C}^2$ form different distributive lattices $(\mathcal{L}^{(1)}, \le)$, $(\mathcal{L}^{(2)}, \le)$ and $(\mathcal{L}^{(3)}, \le)$ whose nontrivial elements, i.e., $\mathrm{ran}(\hat{P}^{(W)}_n) \in \mathcal{L}^{(W)}$ and $\mathrm{ran}(\hat{P}^{(W^\prime \neq W)}_m) \in \mathcal{L}^{(W^\prime \neq W)}$, associated with the different contexts – i.e., collections of all the projection operators which are orthogonal to each other and resolve to the identity operator, cannot meet each other. This corresponds to textbook quantum theory saying that the propositions represented by the incommutable projection operators $\hat{P}^{(W)}_n$ and $\hat{P}^{(W^\prime \neq W)}_m$ cannot be simultaneously verified. Due to such \textit{a contextuality of the Hilbert space}, there does not exist a bivaluation map on the different lattices, i.e., $v: \left\{ \mathcal{L}^{(1)},\mathcal{L}^{(2)},\mathcal{L}^{(3)}  \right\}\;\;\;\cancel{\!\!\mapsto\!\!}\;\;\;\{0,1\}$, even in the two-dimensional case.\\

\bibliographystyle{References}
\bibliography{Hilbert_Lattice_Ref}

\end{document}